\documentclass[a4paper,USenglish,cleveref,autoref,thm-restate]{lipics-v2021}
\hideLIPIcs\def\subjclassHeading{}\ccsdesc{}\global\renewcommand\ccsdesc[2][100]{}
\nolinenumbers 
\title{Mrs.\ Correct and Majority Colorings}

\author{Marcin Anholcer}
{Institute of Informatics and Electronic Economics, Pozna\'n University of Economics and Business, Pozna\'n, Poland}
{m.anholcer@ue.poznan.pl}
{https://orcid.org/0000-0001-7322-7095}
{Partially supported by the National Science Centre of Poland under grant no.\ 2020/37/B/ST1/03298.}

\author{Bart{\l}omiej Bosek}
{Institute of Theoretical Computer Science, Faculty of Mathematics and Computer Science, Jagiellonian University, Krak{\'o}w, Poland}
{bartlomiej.bosek@uj.edu.pl}
{https://orcid.org/0000-0001-8756-3663}
{Partially supported by the National Science Centre of Poland under grant no.\ 2020/37/B/ST1/03298.}

\author{Jaros{\l}aw Grytczuk}
{Faculty of Mathematics and Information Science, Warsaw University of Technology, Warsaw, Poland}
{j.grytczuk@mini.pw.edu.pl}
{https://orcid.org/0000-0002-0258-6143}
{Partially supported by the National Science Centre of Poland under grant no.\ 2020/37/B/ST1/03298.}

\author{Grzegorz Gutowski}
{Institute of Theoretical Computer Science, Faculty of Mathematics and Computer Science, Jagiellonian University, Krak{\'o}w, Poland}
{grzegorz.gutowski@uj.edu.pl}
{https://orcid.org/0000-0003-3313-1237}
{Partially supported by the National Science Centre of Poland under grant no.\ 2019/35/B/ST6/02472.}

\author{Jakub Przyby{\l}o}
{AGH University of Science and Technology, Faculty of Applied Mathematics, al.\,A.\,Mickiewicza\,30, 30-059 Krak{\'o}w, Poland}
{jakubprz@agh.edu.pl}
{https://orcid.org/0000-0002-1262-7017}
{}

\author{Mariusz Zaj\k{a}c}
{Faculty of Mathematics and Information Science, Warsaw University of Technology, Warsaw, Poland}
{m.zajac@mini.pw.edu.pl}
{https://orcid.org/0000-0002-2080-9523}
{Partially supported by the National Science Centre of Poland under grant no.\ 2019/35/B/ST6/02472.}

\authorrunning{M. Anholcer, B. Bosek, J. Grytczuk, G. Gutowski, J. Przyby{\l}o, M. Zaj\k{a}c}
\Copyright{Marcin Anholcer, Bart{\l}omiej Bosek, Jaros{\l}aw Grytczuk, Grzegorz Gutowski, Jakub Przyby{\l}o, Mariusz Zaj\k{a}c}


\usepackage{tikz}
\usetikzlibrary{arrows,math,patterns,shapes}
\usepackage{todonotes}

\let\leq\leqslant
\let\geq\geqslant

\let\rho\varrho

\newcommand{\brac}[1]{{\left(#1\right)}}

\newcommand{\set}[1]{\left\{#1\right\}}

\newcommand{\norm}[1]{{\left|#1\right|}}

\newcommand{\ceil}[1]{{\left\lceil #1 \right\rceil}}

\DeclareMathOperator{\Xtotal}{total}
\DeclareMathOperator{\Xcost}{cost}
\DeclareMathOperator{\Yleft}{left}
\DeclareMathOperator{\Yright}{right}

\begin{document}

\maketitle

\begin{abstract}

A \emph{majority coloring} of a directed graph is a vertex coloring in which each vertex has the same color as at most half of its out-neighbors.
In this note we simplify some proof techniques and generalize previously known results on various generalizations of majority coloring.
In particular, our unified and simplified approach works for \emph{paintability} -- an on-line analog of the list coloring.
\vskip\topmattervskip\baselineskip
\noindent\textcolor{lipicsGray}{\fontsize{9}{12}\sffamily\bfseries MSC2020\enskip} 05C15
\vskip-\topmattervskip\baselineskip
\end{abstract}

\section{Introduction}
\label{sec:intro}

Let $D$ be a finite, directed, simple, loopless graph.
Let $d^+(v)$ denote the number of out-neighbors of a vertex $v$.
A coloring $c$ of the vertices of $D$ is called a \emph{majority coloring} if for every vertex $v \in V$, the number of its out-neighbors in color $c(v)$ is at most $\frac{1}{2}d^+(v)$.
This concept was studied by Kreutzer, Oum, Seymour, van der Zypen, and Wood~\cite{KreutzerOSZW2017}.
It is proved there, among other results, that every directed graph is majority 4-colorable.
It is conjectured that every directed graph is majority 3-colorable and this would be the best possible.

The simple concept of majority coloring can be generalized in various ways.
In the next few subsections we describe different generalizations that are important in our context.

\subsection{Coloring with tolerance}

A \emph{positively-edge-weighted directed graph} $D$ is a tuple $(V,E,\omega)$ where $V=V(D)$ is the set of vertices, $E=E(D)$ is the set of directed edges, each of which is an ordered pair of two different vertices, and $\omega$ is a function that assigns a positive real weight $\omega(vw)$ to every directed edge $vw \in E$.
We can further generalize the majority coloring to positively-edge-weighted directed graphs.
We can also consider a generalization where different vertices allow different proportions of its outgoing edges to be monochromatic. 
Let $\tau$ be a function that assigns a real \emph{tolerance} $\tau(v)$ to every vertex $v \in V$.
A coloring $c$ of the vertices of $D$ is called a \emph{$\tau$-majority coloring} of $D$ if for every vertex $v\in V$, we have that at most $\tau(v)$-fraction of $\omega$-weighted outgoing edges of $v$ are monochromatic.
More formally, for every vertex $v\in V$ we have
$$
\sum_{vw\in E, c(v) = c(w)} \omega(vw) \leq \tau(v) \cdot \sum_{vw\in E} \omega(vw)\textrm{.}
$$

Observe that for a value $\tau(v) \geq 1$ the constraint for vertex $v$ is always satisfied, and that for a value $\tau(v) < 0$ the constraint is never satisfied.
Thus, we are usually interested in the values of tolerance between $0$ and $1$.
When $\tau$ assigns the value $t$ uniformly for every vertex $v \in V$, a $\tau$-majority coloring is called a \emph{$t$-majority coloring}.
Observe that the original majority coloring of a directed graph $(V,E)$ is the $\frac{1}{2}$-majority coloring of the positively edge weighted directed graph $(V,E,\omega)$ where $\omega$ assigns the value $1$ uniformly to every edge $e \in E$.

\subsection{List coloring}

Another generalization is to consider list colorings, where each vertex is colored with a color from a prescribed list.
Suppose that each vertex $v$ is assigned with a list of colors $L(v)$.
We say that $D$ is \emph{$\tau$-majority colorable from $L$} if there is a $\tau$-majority coloring $c$ of $D$ with $c(v) \in L(v)$ for every vertex $v \in V$.
If $D$ is $\tau$-majority colorable from any lists that satisfy $\norm{L(v)} = k$ for every $v \in V$, then we say that $D$ is \emph{$\tau$-majority $k$-choosable}.
It is easy to see that if a graph is $\tau$-majority $k$-choosable, then it is $\tau$-majority $k$-colorable.
Anholcer, Bosek and Grytczuk~\cite{AnholcerBG2017} showed that every directed graph is $\frac{1}{2}$-majority $4$-choosable.
Their technique can also give that for every integer $k \geq 1$, every directed graph is $\frac{1}{k}$-majority $k^2$-choosable.
Later, Gir{\~a}o, Kittipassorn and Popielarz~\cite{GiraoKP2017}, and independently Knox and {\v{S}}{\'a}mal~\cite{KnoxS2018}, showed that every directed graph is $\frac{1}{k}$-majority $2k$-choosable.
It is possible that every directed graph is $\frac{1}{k}$-majority $(2k-1)$-choosable.
This would be the best possible and it is unknown if it holds even for $k=2$.
Some evidence supporting this conjecture was given by Anastos, Lamaison, Steiner and Szab{\'o}~\cite{AnastosLSS2021}.

\subsection{Paintability}

A natural generalization of list colorings is the concept of \emph{paintability} (also called \emph{on-line list coloring}, see~\cite{Schauz2009,Zhu2009}).
Let $\kappa$ be a function that assigns a positive integer $\kappa(v)$ to every vertex $v$ of a positively-edge-weighted directed graph $D=(V,E,\omega)$.
The \emph{$\tau$-majority $\kappa$-painting game on $D$} is a game played in rounds by two players: Lister and Painter (Mr. Paint and Mrs. Correct in the original article of Schauz \cite{Schauz2009}). 
The game starts with all the vertices of $D$ being uncolored.
The $i$-th round, for $i=1,2,\ldots$, starts with Lister presenting a non-empty subset $X_i$ of uncolored vertices of $D$.
Then, Painter responds by selecting a subset $Y_i \subseteq X_i$ and assigns color $i$ to all the vertices in $Y_i$.
Coloring of the vertices in $Y_i$ must satisfy the $\tau$-majority constraints, i.e., the inequality
$$
 \sum_{vw \in E, w \in Y_i} \omega(vw) \leq \tau(v) \cdot \sum_{vw \in E} \omega(vw)
    \textrm{}
$$
must hold for every vertex $v\in Y_i$.
Painter's goal is to color all vertices and Lister's goal is to get some vertex to be uncolored despite multiple presentations.
Specifically, Painter wins the game after the $m$-th round if all the vertices are colored, i.e., $\bigcup_{i=1}^m Y_i = V$.
Lister wins the game after the $m$-th round if there is a vertex $v$ that was presented $\kappa(v)$ many times, and remains uncolored, i.e., for some vertex $v\in V$,
$$
\left|\set{i:v \in X_i}\right| = \kappa(v)
\textrm{ and }
v \notin \bigcup_{i=1}^m Y_i
    \textrm{.}
$$
If neither player wins, the game continues to the next round.
Observe that the requirement that each $X_i$ contains at least one uncolored vertex means that Painter can color at least one vertex in every round, and as a consequence we get that one of the players has a winning strategy in the game.
We say that $D$ is \emph{$\tau$-majority $\kappa$-paintable} if Painter has a winning strategy in the corresponding painting game on $D$.
When $\kappa$ assigns the value $k$ uniformly to every vertex $v \in V$, we say that $D$ is \emph{$\tau$-majority $k$-paintable}.

The idea of paintability generalizes the idea of choosability in the following way.
For any list assignment $L$ that uses all the colors $1,2,\ldots,m$, we can consider a strategy $\mathcal{S}_L$ for Lister that in the $i$-th round presents the set
$$X_i = \set{v \in V: i \in L(v)\text{ and $v$ is uncolored before the $i$-th round}}.$$ (If $X_i$ happens to be empty, then Lister skips presenting $X_i$ and presents the next set).
For $\kappa(v) = \norm{L(v)}$, the winning strategy for Painter in $\tau$-majority $\kappa$-painting game on $D$ played against $\mathcal{S}_L$ yields a $\tau$-majority coloring of $D$ from $L$.
Thus, when $\kappa(v)$ assigns the value $k$ uniformly for every vertex $v \in V$, we get that every $\tau$-majority $k$-paintable graph is $\tau$-majority $k$-choosable.

\subsection{Color ranks}

There is another generalization that we are going to consider in this paper.
It was initially introduced by Anholcer, Bosek and Grytczuk~\cite{AnholcerBG2017} and called \emph{colors with ranks}.
In this variant, each vertex $v \in V$ has a different upper bound on the number of monochromatic out-edges depending on the color assigned to $v$.
More specifically, each color $i$ in the list $L(v)$ has an arbitrarily assigned real number $r_i(v)$, called the \emph{rank} of $i$ in $L(v)$.
Now, the majority-like coloring $c$ with respect to these ranks demands that for each $v$, the number of monochromatic outgoing edges is bounded from above by the rank $r_{c(v)}(v)$.
It was proved in \cite{AnholcerBG2017} that the desired coloring exists whenever the sum of ranks of colors in each list $L(v)$ is at least $2d^+(v)$, which clearly implies the $4$-choosability result mentioned above.

In our setting, it is more convenient to use an analog of the tolerance function to define the same concept of the color rank.
To be more specific, for every possible color $i$, and every vertex $v\in V$, given the \emph{$i$-tolerance of $v$}, $\tau_i(v)$, the respective majority constraints are
$$
    \sum_{vw \in E, c(v)=c(w)} \omega(vw) \leq \tau_{c(v)}(v) \cdot \sum_{vw \in E} \omega(vw)
    \textrm{.}
$$

\subsection{Paintability with ranks}

Finally, we can define a new concept that generalizes both the idea of paintability, and the idea of colors with ranks.
The \emph{ranked-majority $\lambda$-painting game on $D$} is played similarly to the previous painting game.
The $i$-th round starts with Lister presenting a non-empty subset $X_i$ of uncolored vertices of $D$, and a tolerance function $\tau_i$ that assigns a non-negative real value $\tau_i(v)$ to every vertex $v \in X_i$.
Then, Painter responds by selecting a subset $Y_i \subseteq X_i$ and assigns color $i$ to all the vertices in $Y_i$.
The vertices in $Y_i$ must satisfy the majority constraints, i.e.,
$$
    \sum_{vw \in E, w \in Y_i} \omega(vw) \leq \tau_i(v) \cdot \sum_{vw \in E} \omega(vw)
    \textrm{}
$$
must hold for all $v \in Y_i$.
Painter wins the game after the $m$-th round if all the vertices are colored, i.e., $\bigcup_{i=1}^m Y_i = V$.
Lister wins the game after the $m$-th round if there is a vertex $v$ that was presented in sets with total tolerance at least $\lambda(v)$ and remains uncolored, i.e.,
$$
\sum_{i:v \in X_i} \tau_i(v) \geq \lambda(v)
\textrm{ and }
v \notin \bigcup_{i=1,\ldots,m} Y_i
    \textrm{,}
$$
for some vertex $v\in V$.

Observe that the requirement that each $X_i$ contains at least one uncolored vertex with non-negative tolerance means that Painter can color at least one vertex in every round, and as a consequence we get that one of the players has a winning strategy in the game.
Similarly as before, we say that $D$ is \emph{ranked-majority $\lambda$-paintable} if Painter has a winning strategy in the corresponding game.

The ranked-majority paintability generalizes the previously defined majority paintability in the following way.
When Lister presents $\tau_i(v)$ equal to a predefined value $\tau(v)$ independently of $i$, and $\lambda(v) = \kappa(v) \cdot \tau(v)$, then the ranked-majority $\lambda$-painting game on $D$ is the same as the $\tau$-majority $\kappa$-painting game on $D$.

The definition of the ranked-majority painting game requires Lister to present in each round a non-empty subset of uncolored vertices, each with a non-negative tolerance.
It is convenient to relax this requirement and design Lister strategies that can present colored vertices, vertices with negative tolerance, or an empty set of vertices.
We allow for that under the condition that it can happen only finitely many times.
Then, we can modify such a Lister strategy so that colored vertices and vertices with negative tolerance are removed from each presented set.
Further, when Lister strategy presents an empty set of vertices, we can skip it and instead present the set that would be presented in the next round.
The modified strategy is a valid strategy for Lister in the ranked-majority painting game.

\subsection{Undirected graphs}

We have defined the majority coloring concepts for directed graphs, however they were initially considered for undirected graphs, where the majority condition states that every vertex needs a fraction of its neighbors to have a different color than its own.
Lov\'asz~\cite{Lovasz1966} proved that, for every integer $k\geqslant 2$, every finite undirected graph is $\frac{1}{k}$-majority $k$-colorable.
This proof is very simple (take a coloring with the least number of monochromatic edges) and extends easily to the list version giving that any undirected graph is $\frac{1}{k}$-majority $k$-choosable.
A ranked version (for $k=2$) was obtained by Bernardi \cite{Bernardi1987}.
Notice that a \emph{positively-edge-weighted undirected graph} is the same as a symmetric positively-edge-weighted directed graph in which $\omega(vw)=\omega(wv)$ for every edge $vw\in E$.
So, our setting covers all undirected analogs of the introduced variants of the majority coloring.

\subsection{The main results}

The main results of this paper are the following positive results on ranked-majority paintability of undirected, and directed graphs.
\begin{restatable*}[Undirected, Ranked Paintability]{theorem}{thmundirected}\label{thm:undirected}
    Every positively-edge-weighted undirected graph is ranked-majority $1$-paintable.
\end{restatable*}
\begin{restatable*}[Directed, Ranked Paintability]{theorem}{thmdirected}\label{thm:directed}
    Every positively-edge-weighted directed graph is ranked-majority $2$-paintable.
\end{restatable*}

The proofs of these results consist of three main parts. The first one (\autoref{lem:main}) develops the ideas of Lov\'asz~\cite{Lovasz1966}. The second one is an application of the celebrated Perron-Frobenius Theorem---a powerful tool used in various disciplines, like Computer Science (search engines), Demography (population growth models), Economy (growth and equilibrium models) or Statistics (Markov chains, random walks). We use it mainly to transfer the results from undirected graphs to strongly connected directed graphs (\autoref{lem:scc}). This idea was applied in the context of majority colorings by Knox and {\v{S}}{\'a}mal~\cite{KnoxS2018}, and by Gir{\~a}o, Kittipassorn and Popielarz~\cite{GiraoKP2017}, who used a special case of the result formulated previously by Alon~\cite{AlonSplitting}. Finally, the third part (\autoref{lem:ranks}) develops the idea of ranked colors from \cite{AnholcerBG2017} which allows to extend the result to all directed graphs.

It is worth noticing that these two theorems generalize and improve on the previously known results on majority colorings in both the list coloring setting, and the colors with ranks setting.
These previously known results follow as easy corollaries which are presented in \autoref{sec:un}, and \autoref{sec:dir}. Simultaneously, the more general statement of the problems makes some parts of the proofs easier. We believe that the proofs presented in this paper are easier to follow than the previous ones.


\section{Results and proofs}

\subsection{Undirected graphs}\label{sec:un}

First we prove the following lemma that allows for a construction of easy, yet effective, strategies for Painter in majority painting games played on positively-edge-weighted undirected graphs.
Intuitively, the lemma gives a good Painter response $Y$ to any Lister move $X$.
For every vertex $v$ in $X$ we have that: either $v$ is in $Y$ (and gets colored), or $v$ cannot be added to $Y$ because too many neighbors of $v$ are in $Y$ (and many neighbors of $v$ get colored).
The proof is based on a similar idea as the argument of Lov\'asz~\cite{Lovasz1966}: we choose the set maximizing a certain \emph{potential} function and then prove that moving a vertex into or out of this set would lead to a contradiction.

\begin{lemma}\label{lem:main}
Let $G=(V,E,\omega)$ be a positively-edge-weighted undirected graph.
Let $X\subseteq V$ be a subset of vertices and $\rho$ be a vertex ranking that assigns a real rank $\rho(v)$ to every vertex $v \in X$.
Then, there exists a subset $Y \subseteq X$ such that for every $v \in X$ we have
$$
    v \in Y \iff \brac{\rho(v) \geq \sum_{vw \in E, w \in Y} \omega(vw)}
    \textrm{.}
$$
\end{lemma}

\begin{proof}
Given the set $X$, for any vertex $v \in X$ we define the cost of $v$ in the following way:
$$
\begin{aligned}
    \Xtotal_v & = \sum_{vw\in E,w\in X} \omega(vw)\textrm{,}\\
    \Xcost_v & = 2\rho(v) - \Xtotal_v\textrm{.}
\end{aligned}
$$
Now, given any subset $Y \subseteq X$ we can define the cost of $Y$ as the sum of the costs of vertices in $Y$ plus the sum of weights of all the edges in the cut $(Y,X\setminus{}Y)$:
$$
    \Xcost_Y = \sum_{v \in Y} \Xcost_v + \sum_{vw\in E, v \in Y, w \in X\setminus{}Y} \omega(vw)\textrm{.}
$$

Now, let $Y$ be such that it maximizes the value of cost among all subsets of $X$ and maximizes the size of $Y$ among all subsets of $X$ that maximize the value of cost.
We shall prove that $Y$ satisfies the conditions of the lemma.
Suppose to the contrary that for some vertex $v \in X$ the condition does not hold.
Define:
$$
\begin{aligned}
\Yleft_v & = \sum_{vw \in E, w \in Y} \omega(vw)\textrm{,}\\
\Yright_v & = \sum_{vw \in E,w \in X\setminus{}Y} \omega(vw)\textrm{,}
\end{aligned}
$$
and observe that $\Yleft_v + \Yright_v = \Xtotal_v$.
We distinguish two cases, depending on whether $v \in Y$ or not.

For the first case, suppose that $v \in Y$ and $\rho(v) < \Yleft_v$.
Consider the set $Y^\prime = Y\setminus{}\set{v}$ and observe that: 
$$
\begin{aligned}
    \Xcost_{Y^\prime} & = \Xcost_Y - \Xcost_v - \Yright_v + \Yleft_v  \\
    & = \Xcost_Y - (2\rho(v) - \Xtotal_v) - (\Xtotal_v-\Yleft_v) + \Yleft_v  \\
    & = \Xcost_Y - 2\rho(v) + 2\Yleft_v  \\
    & > \Xcost_Y - 2\rho(v) + 2\rho(v)  \\
    & = \Xcost_Y\textrm{,}
\end{aligned}
$$ 
which contradicts the choice of maximal $Y$, as we have $\Xcost_{Y^\prime} > \Xcost_{Y}$.

Similarly for the second case, suppose that $v \notin Y$ and $\rho(v) \geq \Yleft_v$.
Consider the set $Y^\prime = Y\cup\set{v}$ and observe that:
$$
\begin{aligned}
    \Xcost_{Y^\prime} & = \Xcost_Y + \Xcost_v + \Yright_v - \Yleft_v  \\
    & = \Xcost_Y + (2\rho(v) - \Xtotal_v) + (\Xtotal_v-\Yleft_v) - \Yleft_v  \\
    & = \Xcost_Y + 2\rho(v) - 2\Yleft_v  \\
    & \geq \Xcost_Y + 2\rho(v) - 2\rho(v)  \\
    & = \Xcost_Y\textrm{,}
\end{aligned}
$$
which contradicts the choice of maximal $Y$, as we have $\Xcost_{Y^\prime} \geq \Xcost_{Y}$ and $\norm{Y^\prime} > \norm{Y}$.
\end{proof}

Applying \autoref{lem:main} directly in every round of the painting game gives us the following result.
\thmundirected

\begin{proof}
Let $G=(V,E,\omega)$ be a positively-edge-weighted undirected graph.
The strategy for the Painter in the ranked-majority $1$-painting game on $G$ is as follows.
In the $i$-th round, when Lister presents $X_i$ and $\tau_i$, Painter defines $\rho_i(v)=\tau_i(v) \cdot \sum_{vw \in E} \omega(vw)$, applies \autoref{lem:main} for $X=X_i$ and $\rho=\rho_i$, and obtains set $Y_i = Y$.
This lemma guarantees that $Y_i$ satisfies the majority constraints, since for each $v\in Y_i$,
$$
\sum_{vw \in E, w \in Y_i} \omega(vw) \leq \rho_i(v) = \tau_i(v) \cdot \sum_{vw \in E} \omega(vw)
\textrm{.}
$$

Now, assume for a contradiction that Lister wins the game after the $m$-th round and for some vertex $v$ we have
$$
\sum_{i:v \in X_i} \tau_i(v) \geq 1
\textrm{ and }
v \notin \bigcup_{i=1}^m Y_i
    \textrm{.}
$$
As Painter did not color vertex $v$, we get
$$
\sum_{vw \in E, w \in Y_i} \omega(vw) > \rho_i(v) = \tau_i(v) \cdot \sum_{vw \in E} \omega(vw)
    \textrm{,}
$$
for every color $i$ such that $v\in X_i$.
Summing over all $i$ with $v \in X_i$ we get:
$$
 \sum_{vw \in E}  \omega(vw) \geq \sum_{i:v \in X_i}\brac{\sum_{vw \in E, w \in Y_i} \omega(vw)} > \sum_{i:v \in X_i} \brac{\tau_i(v) \cdot \sum_{vw \in E} \omega(vw)} \geq 1 \cdot \sum_{vw \in E}  \omega(vw)
    \textrm{,}
$$
which gives a contradiction.
Thus, when the game ends, every vertex is colored and Painter wins.
\end{proof}

The following corollary of \autoref{thm:undirected} easily generalizes the result by Anholcer, Bosek and Grytczuk~\cite{AnholcerBG2017} on choosability with ranked colors.
\begin{corollary}[Undirected, Ranked Choosability]\label{cor:un_rank_choice}
    Let $G=(V,E,\omega)$ be a positively-edge-weighted undirected graph.
    Suppose that each vertex $v$ is assigned a list $L(v)$ of colors.
    Suppose that for each vertex $v$, each color $i\in L(v)$ is assigned a real number $r_i(v)$, the \emph{rank} of color $i$ in $L(v)$.
    Assume further that for every vertex $v$, the color ranks $r_i(v)$ satisfy the following condition:
    $$
    \sum_{i \in L(v)} r_i(v) \geq \sum_{vw\in E} \omega(vw)\textrm{.}
    $$
    Then there is a vertex coloring of $G$ from lists $L(v)$ satisfying the following constraint: If $v$ is colored by $i$, then the sum of weights of the edges connecting $v$ to its neighbors in the same color $i$ is at most $r_i(v)$.
\end{corollary}
\begin{proof}
    Assume, without loss of generality, that 
    $$
    \bigcup_{v\in V} L(v) = \{1,\ldots, m\}
    $$
    for some positive integer $m$.
    Consider the following strategy $\mathcal{S}$ for Lister in ranked-majority $1$-painting game on $G$.
    In the $i$-th round of the game Lister presents set the $X_i = \set{v \in V: i \in L(v)}$, and the $i$-tolerance function given by
    $$
    \tau_i(v) = \frac{r_i(v)}{\sum_{vw\in E} \omega(vw)}
    \text{.}
    $$
    The winning strategy for Painter provided by \autoref{thm:undirected} played against $\mathcal{S}$ constructs the desired coloring of $G$.
\end{proof}

Simplifying \autoref{thm:undirected}, when the tolerance function does not depend on $i$, we get the following.

\begin{corollary}[Undirected, Non-uniform Paintability]\label{cor:un_nonu_paint}
    Let $G$ be a positively-edge-weighted undirected graph.
    Let $\tau$ be a tolerance function that assigns $0 \leq \tau(v) \leq 1$ for every vertex $v$ of $G$.
    For $\kappa(v) = \ceil{\frac{1}{\tau(v)}}$, $G$ is $\tau$-majority $\kappa$-paintable.
\end{corollary}

Simplifying further, when $\tau$ assigns the value $\frac{1}{k}$ uniformly for every vertex $v \in V$, we get the following new generalization of the previously known results.

\begin{corollary}[Undirected Paintability]\label{cor:un_paint}
    Every positively-edge-weighted undirected graph is $\frac{1}{k}$-majority $k$-paintable.
\end{corollary}

\subsection{Strongly connected directed graphs}\label{sec:scc}

Now, we move from undirected graphs to directed graphs.
First, we show how to deal with strongly connected directed graphs.
Our approach here follows the ideas from the works of Gir{\~a}o, Kittipassorn and Popielarz~\cite{GiraoKP2017}, and of Knox and {\v{S}}{\'a}mal~\cite{KnoxS2018}, and uses Perron-Frobenius Theorem to transform a strongly connected directed graph into an undirected one.
This transformation comes with a price of a multiplicative factor $2$, as described in the following lemma.

\begin{lemma}\label{lem:scc}
    Every strongly connected positively-edge-weighted directed graph is ranked-majority $2$-paintable.
\end{lemma}

\begin{proof}
Let $D=(V,E_D,\omega_D)$ be a strongly connected positively-edge-weighted directed graph.
Observe that, since in the directed case the incoming edges do not influence the majority condition, we can modify original weights $\omega_D$ by uniformly scaling weights of all outgoing edges of a single vertex $v \in V$ in the following way:
$$
\omega_D(vw):=\frac{\omega_D(vw)}{\sum_{u:vu\in E}{\omega_D(vu)}}\text{,}
$$
and get the same $\tau$-majority colorings of the modified graph.
Thus, we can assume that for every vertex the total weight of its outgoing edges equals $1$.
We can assume that $V=\set{1,2,\ldots,n}$, and define $T$ to be an $n\times n$ matrix so that $T_{vw}=\omega(vw)$ for every edge $vw \in E$, and $T_{vw}=0$ for all other entries.
We have that $T_{vw}$ is non-negative, and for every $v\in V$, we have
$$
\sum_{w\in V}T_{vw}=1.
$$

When $Ty=cy$ for some real $c$ and a real vector $y$, then by choosing $v$ with the maximum value of $|y_v|$ we get
$$
|(cy)_v|=|\sum_{w\in V}T_{vw}y_w|\leq
\sum_{w\in V}|T_{vw}||y_w|\leq
\sum_{w\in V}T_{vw}|y_v|=|y_v|,
$$
so $|c|\leq 1$, which means that the spectral radius (i.e., the maximum of the absolute values of the eigenvalues) of $T$ is $1$.

Recall that a square matrix is irreducible if it cannot be transformed into a block upper triangular form by permuting its rows and columns. In particular, if one replaces all the non-zero entries of a matrix with $1$'s and considers it as an adjacency matrix of a digraph, it is irreducible if and only if this digraph is strongly connected (see e.g.~\cite[Section 8.7]{refGodsil}). Since $D$ is strongly connected, $T$ is irreducible and we can apply the Perron-Frobenius Theorem (see e.g.~\cite[Theorem 8.8.1]{refGodsil}).
According to it, there exists a positive left eigenvetor $x$ of $T$ with eigenvalue $1$, that is a vector $x$ such that $x_v>0$ for every $v\in V$, and
$$
\sum_{w\in V}T_{wv}x_w=x_v.
$$

Consider now the underlying positively-edge-weighted undirected graph $G=(V,E_G,\omega_G)$ of $D$, with weighting $\omega_G$ defined as follows:
$$
\omega_G(vw) = \omega_G(wv) = x_vT_{vw} + x_wT_{wv}
\textrm{.}
$$
Observe that the total weight of edges incident to any vertex $v$ equals
$$
\sum_{vw\in E_G} \omega_G(vw) = \sum_{vw\in E_G} (x_vT_{vw} + x_wT_{wv}) = \sum_{w \in V} x_vT_{vw} + \sum_{w \in V} x_wT_{wv} = x_v + x_v = 2x_v
\textrm{.}
$$
Let $\mathcal{S}$ be the winning strategy for Painter in the ranked-majority $1$-painting game on $G$ guaranteed by \autoref{thm:undirected}.

Now, we define the winning strategy for Painter in the \emph{real game} -- the ranked-majority $2$-painting game on $D$.
In order to use strategy $\mathcal{S}$ as a \emph{subprocedure}, our new strategy plays as Lister a \emph{side game} of ranked-majority $1$-painting game on $G$ against $\mathcal{S}$.
Lister is bound to lose this side game, and the key idea is to encode the real game in the side game in order to beat Lister in the real game.

In the $i$-th round of the real game, when Lister presents the set $X_i$, and the $i$-tolerance function $\tau_i$, Painter constructs a function $\tau^\prime_i(v) = \frac{\tau_i(v)}{2}$,
plays the move $X_i$, $\tau^\prime_i$ as Lister against strategy $\mathcal{S}$ in the side game, and copies the response $Y_i$ as a response in the real game.

The set $Y_i$ is a valid Painter response in the side game, which, for every $v\in Y_i$, satisfies the following inequality:
$$
    \sum_{vw \in E_G, w \in Y_i} \omega_G(vw) \leq \tau^\prime_i(v) \cdot \sum_{vw \in E_G, w \in V} \omega_G(vw) = \frac{\tau_i(v)}{2} \cdot 2x_v = \tau_i(v) \cdot x_v
    \textrm{.}
$$
Thus, for every $v\in Y_i$,
$$
    \sum_{vw \in E_D, w \in Y_i} \omega_D(vw) \leq \sum_{vw \in E_G, w \in Y_i} \frac{x_vT_{vw}}{x_v} \leq \frac{\sum_{vw \in E_G, w \in Y_i} \omega_G(vw)}{x_v} \leq \tau_i(v)
    \textrm{,}
$$
and hence $Y_i$ is a valid Painter response in the real game.

Now, assume to the contrary that Lister wins the real game, and let $v$ be a vertex that remains uncolored that was presented in sets with total tolerance at least $2$.
The same vertex $v$ was presented in sets with total tolerance at least $1$ and remains uncolored in the side game, which contradicts $\mathcal{S}$ being a winning strategy.
\end{proof}

\subsection{Directed graphs}\label{sec:dir}

When a directed graph is not strongly connected, we can divide it into strongly connected components.
We say that a strongly connected component is a \emph{source} component if there are no edges directed to the component from other components of the graph.
The following technical lemma allows Painter to play the painting game in a source component almost independently of the game played in the rest of the graph.

For a positively-edge-weighted directed graph $D=(V,E,\omega)$, and any subset $W \subseteq V$ of vertices we define $D[W]=(W,E|_W,\omega|_{E|_W})$ to be the positively-edge-weighted directed graph on the vertex set $W$ with both edges, and weight function restricted to $W$.
\begin{lemma}[Rank Reduction]\label{lem:ranks}
   Let $D=(V,E,\omega)$ be a positively-edge-weighted directed graph.
   Let $S$ be the set of vertices of some source strongly connected component of $D$, and let $T=V\setminus S$.
   Let $\lambda$ be a function that assigns a positive real number to every vertex of $D$, with $\lambda(v) \geq 1$ for every $v \in S$.
   If $D[S]$ is ranked-majority $\lambda|_S$-paintable, and $D[T]$ is ranked-majority $\lambda|_T$-paintable, then
   $D$ is ranked-majority $\lambda$-paintable.
\end{lemma}

\begin{proof}
We construct a winning strategy for Painter in the \emph{real game} -- ranked-majority $\lambda$-painting game on $D$, as follows.
There are two side games.
A $T$-game, where we play as Lister a ranked-majority $\lambda|_T$-painting game on $D[T]$ against a winning strategy $\mathcal{S}_T$.
And an $S$-game, where we play as Lister a ranked-majority $\lambda|_S$-painting game on $D[S]$ against a winning strategy $\mathcal{S}_S$.

In the $i$-th round of the real game, when Lister presents the set $X_i$, and the $i$-tolerance function $\tau_i$,
Painter first plays the move $X_i \cap T$, $\tau_i$ as Lister in the $T$-game and obtains the response $Y_{i,T}$.
For every vertex $v \in S$, let 
$$\rho_i(v) = \tau_i(v)\cdot\sum_{vw \in E}\omega(vw) - \sum_{vw \in E,w \in Y_{i,T}}\omega(vw)$$
denote the maximum total weight of monochromatic edges from $v$ to other vertices in $S$ assuming that $Y_i\cap T = Y_{i,T}$.
For every vertex $v \in S$, let 
$$\tau^\prime_i(v) = \frac{\rho_i(v)}{\sum_{vw \in E, w \in S}\omega(vw)}\textrm{.}$$
Painter plays the move $X_i \cap S$, $\tau^\prime_i$ as Lister in the $S$-game and obtains the response $Y_{i,S}$.
Painter responds with the set $Y_i = Y_{i,T} \cup Y_{i,S}$ in the real game.

Since there are no edges from $T$ to $S$, and by the definition of $\tau^\prime_i$ it is easy to see that $Y_i$ satisfies the necessary majority constraints.
Now, assume to the contrary that Lister wins the game, and let $v$ be a vertex that was presented in the sets with total tolerance at least $\lambda(v)$ and remains uncolored.
It is easy to see that $v \notin T$ is guaranteed by the fact that $\mathcal{S}_T$ is a winning strategy.
Thus, $v \in S$, and $\sum_{i:v\in X_i}\tau_i(v) \geq \lambda(v)$.
So, using the inequality $\lambda(v)\geq 1$ we obtain
$$
\sum_{i:v\in X_i} \rho_i(v) \geq \lambda(v)\cdot\sum_{vw\in E}\omega(vw) - \sum_{vw \in E, w \in T} \omega(vw) \geq \lambda(v)\cdot\sum_{vw \in E, w \in S} \omega(vw)
\textrm{,}
$$
and
$$
\sum_{i:v \in X_i} \tau^\prime_i(v) \geq \lambda(v)
$$
gives a contradiction with $\mathcal{S}_S$ being a winning strategy.
\end{proof}

\autoref{lem:ranks} allows us to extend the result from \autoref{lem:scc} to all directed graphs, not necessarily strongly connected.
\thmdirected

\begin{proof}
Let $D=(V,E,\omega)$ be a positively-edge-weighted directed graph.
The proof is by induction on the number of strongly connected components of $D$.
If $D$ is strongly connected, the result follows from \autoref{lem:scc}.
Otherwise, let $S$ be the vertex set of some source strongly connected component of $D$, and $T=V\setminus S$.
$D[S]$ is strongly connected and ranked-majority $2$-paintable by \autoref{lem:scc}.
The number of strongly connected components in $D[T]$ is smaller than in $D$, and $D[T]$ is ranked-majority $2$-paintable by induction.
Thus, $D$ is ranked-majority $2$-paintable by \autoref{lem:ranks}.
\end{proof}

\autoref{thm:directed} immediately implies the following analogs of Corollaries~\ref{cor:un_rank_choice},~\ref{cor:un_nonu_paint} and~\ref{cor:un_paint}.

\begin{corollary}[Directed, Ranked Choosability]\label{cor:dir_rank_choice}
    Let $D$ be a positively-edge-weighted directed graph.
    Suppose that each vertex $v$ is assigned with a list $L(v)$ of colors.
    Suppose further that for each vertex $v$, for each color $i$ in $L(v)$ there is a real number $r_i(v)$, the \emph{rank} of color $i$ in $L(v)$.
    Assume that for every vertex $v$, the color ranks $r_i(v)$ satisfy the following condition:
    $$
    \sum_{i \in L(v)} r_i(v) \geq 2 \cdot \sum_{vw\in E} \omega(vw)\textrm{.}
    $$
    Then there is a vertex coloring of $G$ from lists $L(v)$ satisfying the following constraint: If $i$ is a color assigned to $v$, then the sum of weights of edges connecting $v$ to an out-neighbor in color $i$ is at most $r_i(v)$.
\end{corollary}

\begin{corollary}[Directed, Non-uniform Paintability]\label{cor:dir_nonu_paint}
    Let $D$ be a positively-edge-weighted directed graph.
    Let $\tau$ be a tolerance function that assigns $0 \leq \tau(v) \leq 1$ for every vertex $v$ of $D$.
    For $\kappa(v) = \ceil{\frac{2}{\tau(v)}}$, $D$ is $\tau$-majority $\kappa$-paintable.
\end{corollary}

\begin{corollary}[Directed Paintability]\label{cor:dir_paint}
    Every positively-edge-weighted directed graph is $\frac{1}{k}$-majority $2k$-paintable.
\end{corollary}


\section{Discussion}

It is easy to see that the complete undirected graph $K_k$ is not $\frac{1}{k}$-majority $(k-1)$-colorable, and is not ranked-majority $\frac{k-1}{k}$-paintable.
Similarly, any orientation of $K_{2k-1}$ with out-degree of every vertex exactly $k-1$ is not $\frac{1}{k}$-majority $(2k-2)$-colorable, and is not ranked-majority $\frac{2k-2}{k}$-paintable.
Thus, the constant $1$ in \autoref{thm:undirected}, and the constant $2$ in \autoref{thm:directed} are optimal.
Still, it is an interesting question for which tolerance functions we can get stronger results.
In particular, it is possible that every positively-edge-weighted directed graph is $\frac{1}{2}$-majority $3$-paintable.

The problems of majority coloring and majority list coloring were also considered for infinite graphs.
The famous \emph{Unfriendly Partition Conjecture} by Cowan and Emerson (\cite{CowEm}, see~\cite{Aharoni1990}) states that every countable undirected graph is majority 2-colorable.
It was proved for graphs with finitely many vertices of infinite degree by Aharoni, Milner and Prikry~\cite{Aharoni1990}, for rayless graphs by Bruhn, Diestel, Georgakopoulos, and Spr\"ussel~\cite{Diestel2010} and for graphs not containing an infinite clique subdivision by Berger~\cite{Berger2017}.
On the other hand, Shelah and Milner~\cite{shelah_milner_1990} showed that every infinite undirected graph is majority 3-colorable and that there are uncountable graphs for which 3 colors are necessary.
Anholcer, Bosek and Grytczuk~\cite{ABGInf2020} proved that every countable directed graph is $\frac{1}{2}$-majority $4$-choosable.
Recently Haslegrave~\cite{Haslegrave2020} showed that every countable undirected graph is $\frac{1}{2}$-majority $3$-choosable and that the same holds for countable directed acyclic graphs. Somewhat surprisingly, the later result turned out to be optimal, as proved recently by Bosek and Katan \cite{BosekKatan2024}. Actually, they constructed an infinite countable directed acyclic graph which is not $\frac{1}{2}$-majority $2$-colorable.

\bibliography{paper}

\begin{thebibliography}{10}

\bibitem{Aharoni1990}
Ron Aharoni, Eric~C. Milner, and Karel Prikry.
\newblock Unfriendly partitions of a graph.
\newblock {\em Journal of Combinatorial Theory Series B}, 50(1):1--10, 1990.
\newblock \href {https://doi.org/10.1016/0095-8956(90)90092-E}
  {\path{doi:10.1016/0095-8956(90)90092-E}}.

\bibitem{AlonSplitting}
Noga Alon.
\newblock Splitting digraphs.
\newblock {\em Combinatorics, Probability and Computing}, 15:933--937, 2006.
\newblock \href {https://doi.org/10.1017/S0963548306008042}
  {\path{doi:10.1017/S0963548306008042}}.

\bibitem{AnastosLSS2021}
Michael Anastos, Ander Lamaison, Raphael Steiner, and Tibor Szab{\'o}.
\newblock Majority colorings of sparse digraphs.
\newblock {\em Electronic Journal of Combinatorics}, 28(2):P2.31:1--17, 2021.
\newblock \href {https://doi.org/10.37236/10067} {\path{doi:10.37236/10067}}.

\bibitem{AnholcerBG2017}
Marcin Anholcer, Bart{\l}omiej Bosek, and Jaros{\l}aw Grytczuk.
\newblock Majority choosability of digraphs.
\newblock {\em Electronic Journal of Combinatorics}, 24(3):P3.57:1--5, 2017.
\newblock \href {https://doi.org/10.37236/6923} {\path{doi:10.37236/6923}}.

\bibitem{ABGInf2020}
Marcin Anholcer, Bart{\l}omiej Bosek, and Jaros{\l}aw Grytczuk.
\newblock Majority choosability of countable graphs.
\newblock {\em European Journal of Combinatorics}, 117:Article 103829: 1--8,
  2024.
\newblock \href {https://doi.org/10.1016/j.ejc.2023.103829}
  {\path{doi:10.1016/j.ejc.2023.103829}}.

\bibitem{Berger2017}
Eli Berger.
\newblock Unfriendly partitions for graphs not containing a subdivision of an
  infinite clique.
\newblock {\em Combinatorica}, 37(2):157--166, 2017.
\newblock \href {https://doi.org/10.1007/s00493-015-3261-1}
  {\path{doi:10.1007/s00493-015-3261-1}}.

\bibitem{Bernardi1987}
Claudio Bernardi.
\newblock On a theorem about vertex colorings of graphs.
\newblock {\em Discrete Mathematics}, 64(3):95--96, 1987.
\newblock \href {https://doi.org/10.1016/0012-365X(87)90243-3}
  {\path{doi:10.1016/0012-365X(87)90243-3}}.

\bibitem{BosekKatan2024}
Bart{\l}omiej Bosek and Aleksander Katan.
\newblock A note about majority colorings of countable dags.
\newblock {\em arXiv}, 2024.
\newblock \href {https://doi.org/10.48550/arXiv.2406.04189}
  {\path{doi:10.48550/arXiv.2406.04189}}.

\bibitem{Diestel2010}
Henning Bruhn, Reinhard Diestel, Agelos Georgakopoulos, and Philipp Spr\"ussel.
\newblock Every rayless graph has an unfriendly partition.
\newblock {\em Combinatorica}, 30(5):521--532, 2010.
\newblock \href {https://doi.org/10.1007/s00493-010-2590-3}
  {\path{doi:10.1007/s00493-010-2590-3}}.

\bibitem{CowEm}
Robert Cowen and William Emerson.
\newblock Proportional colorings of graphs.
\newblock Unpublished.

\bibitem{GiraoKP2017}
Ant{\'o}nio Gir{\~a}o, Teeradej Kittipassorn, and Kamil Popielarz.
\newblock Generalized majority colourings of digraphs.
\newblock {\em Combinatorics, Probability and Computing}, 26(6):850--855, 2017.
\newblock \href {https://doi.org/10.1017/S096354831700044X}
  {\path{doi:10.1017/S096354831700044X}}.

\bibitem{refGodsil}
Chris Godsil and Gordon~F. Royle.
\newblock {\em Algebraic Graph Theory}.
\newblock Graduate Texts in Mathematics. Springer New York, 2013.

\bibitem{Haslegrave2020}
John Haslegrave.
\newblock Countable graphs are majority 3-choosable.
\newblock {\em Discussiones Mathematicae Graph Theory}, 2020.
\newblock \href {https://doi.org/10.7151/dmgt.2383}
  {\path{doi:10.7151/dmgt.2383}}.

\bibitem{KnoxS2018}
Fiachra Knox and Robert {\v{S}}{\'a}mal.
\newblock Linear bound for majority colourings of digraphs.
\newblock {\em Electronic Journal of Combinatorics}, 25(3):P3.29:1--4, 2018.
\newblock \href {https://doi.org/10.37236/6762} {\path{doi:10.37236/6762}}.

\bibitem{KreutzerOSZW2017}
Stephan Kreutzer, Sang-il Oum, Paul Seymour, Dominic van~der Zypen, and
  David~R. Wood.
\newblock Majority colourings of digraphs.
\newblock {\em Electronic Journal of Combinatorics}, 24(2):P2.25:1--9, 2017.
\newblock \href {https://doi.org/10.37236/6410} {\path{doi:10.37236/6410}}.

\bibitem{Lovasz1966}
L\'aszl\'o Lov\'asz.
\newblock On decomposition of graphs.
\newblock {\em Studia Scientiarum Mathematicarum Hungarica}, I(1-2):237--238,
  1966.

\bibitem{Schauz2009}
Uwe Schauz.
\newblock {M}r.~{P}aint and {M}rs.~{C}orrect.
\newblock {\em Electronic Journal of Combinatorics}, 16(1):R77:1--18, 2009.
\newblock \href {https://doi.org/10.37236/166} {\path{doi:10.37236/166}}.

\bibitem{shelah_milner_1990}
Saharon Shelah and Eric~C. Milner.
\newblock Graphs with no unfriendly partitions.
\newblock In {\em A Tribute to Paul Erd{\H{o}}s}, page 373–384. Cambridge
  University Press, 1990.
\newblock \href {https://doi.org/10.1017/CBO9780511983917.031}
  {\path{doi:10.1017/CBO9780511983917.031}}.

\bibitem{Zhu2009}
Xuding Zhu.
\newblock On-line list colouring.
\newblock {\em Electronic Journal of Combinatorics}, 16(1):R127:1--16, 2009.
\newblock \href {https://doi.org/10.37236/216} {\path{doi:10.37236/216}}.

\end{thebibliography}

\end{document}